\documentclass[12pt]{amsart}

\usepackage{epsf,amssymb,times, graphicx}
\newtheorem{thm}{Theorem}[section]

\newtheorem{prop}[thm]{Proposition}

\newtheorem{lemma}[thm]{Lemma}
\newtheorem{cor}[thm]{Corollary}

\newtheorem{question}[thm]{Question}

\newcommand{\R}{\mathbb{R}}
\newcommand{\N}{\mathbb{N}}
\newcommand{\Z}{\mathbb{Z}}

\newcommand{\bdry}{\partial}

\newcommand{\ie}{{\em i.e.}}
\newcommand{\CC}{\mathcal{C}}

\newcommand{\CF}{\mathcal{F}}

\newcommand{\Sg}{\Sigma}
\newcommand{\ra}{\rightarrow}
\newcommand{\hra}{\hookrightarrow}
\newcommand{\co}{\colon\thinspace}

\newcommand{\tb}{tb}
\newcommand{\dv}{div}
\newcommand{\slope}{slope}
\newcommand{\intr}{int}

\newcommand{\s}{\vskip.13in}
\newcommand{\n}{\noindent}

\newcommand{\be}{\begin{enumerate}}
\newcommand{\ee}{\end{enumerate}}
\newcommand{\bp}{\begin{proof}}
\newcommand{\ep}{\end{proof}}

\begin{document}


\title{Contact Structures on Open 3-Manifolds}

\author{James J.\ Tripp}
\address{Department of Mathematics,
University of Pennsylvania,
209 South 33rd St.,
Philadelphia, PA 19104-6395}
\email{jtripp@math.upenn.edu}
\urladdr{http://www.math.upenn.edu/\char126 jtripp}


\date{This Version: August 22, 2004}

\begin{abstract} In this paper, we study contact structures on any open $3$-manifold $V$ which is the interior of a compact
$3$-manifold.  To do this, we introduce proper contact isotopy invariants called the slope at infinity and the division number
at infinity.  We first prove several classification theorems for $T^2 \times [0, \infty)$, $T^2 \times \R$, and $S^1 \times
\R^2$ using these concepts.  This investigation yields infinitely many tight contact structures on $T^2 \times [0,
\infty)$, $T^2 \times \R$, and $S^1 \times \R^2$ which admit no precompact embedding into another tight contact structure
on the same space.  Finally, we show that if $V$ is irreducible and has an end of nonzero genus, then there are uncountably
many tight contact structures on $V$ that are not contactomorphic, yet are isotopic.  Similarly, there
are uncountably many overtwisted contact structures on $V$ that are not contactomorphic, yet are isotopic.\end{abstract}

\maketitle


\section{Introduction}


Recently, there has been much work towards the classification of tight contact structures on compact $3$-manifolds up to
isotopy (relative to the boundary).  In particular, Honda and Giroux provided several classification theorems for solid
tori, toric annuli, torus bundles over the circle, and circle bundles over surfaces \cite{Gi1, Gi2, Gi3, Ho2, Ho3}.  In
comparison, tight contact structures on open $3$-manifolds have been virtually unstudied.  Two main results dealing with
open contact manifolds are due to Eliashberg.  In \cite{El1}, Eliashberg shows that $\R^3$ has a unique tight contact
structure.  It is immediate from his proof that $S^2 \times [0,\infty)$ has a unique tight contact structure with a fixed
characteristic foliation on $S^2 \times 0$.  Therefore, the classification of tight contact structures on open manifolds
with only $S^2$ ends can be reduced to the case of compact manifolds.  In \cite{El3}, Eliashberg shows that, in contrast to
the situation for $S^2$ ends, there are uncountably many tight contact structures on $S^1 \times \R^2$ that are not
contactomorphic.  The situation for closed $3$-manifolds is different.  Colin, Giroux, and Honda proved that an atoroidal $3$-manifold supports finitely many tight contact structures up to isotopy \cite{CGH}.   Honda, Kazez, and Mati\'c, and independently, Colin, show that an irreducible, toroidal $3$-manifold supports countably infinitely many tight contact structures up to isotopy \cite{HKM1, Co}.  

In this paper, we study tight contact structures on any open manifold $V$ which is the interior of a compact manifold.  Due
to the failure of Gray's Theorem on open contact manifolds, we relegate ourselves to the study of tight contact structures
up to {\em proper isotopy}, by which we mean isotopy of the underlying manifold rather than a one-parameter family of
contact structures.  When we say that two contact structures are isotopic, we will mean that they are connected by a
one-parameter family of contact structures.  We first introduce two new proper isotopy invariants which we call the {\em
slope at infinity} and the {\em division number at infinity} of an end $\Sigma_g \times [0, \infty)$ of an open contact
manifold.  These invariants are most naturally defined for toric ends $T^2 \times [0, \infty)$, where we take our
inspiration from the usual definition of the slope and division number of a convex torus.  Using these invariants and
Honda's work in \cite{Ho2}, we essentially classify tight contact structures on toric ends $T^2 \times [0, \infty)$.  In
particular, we show that there is a natural bijection between tight toric annuli and tight toric ends that {\em attain} the
slope at infinity and have finite division number at infinity.  However, we also show that for any slope at infinity there
is an infinite family of tight toric ends which do not attain the slope at infinity and therefore do not come from closed
toric annuli.  Interestingly, these contact structures are strange enough that they cannot be properly embedded in another
tight contact manifold.  This yields the following 

\begin{thm} \label{thm: uncountableembeddings} 

Let $X$ be $T^2 \times [0,1)$, $T^2 \times (0,1)$ or $S^1 \times D^2$, where $D^2$ is the open unit disk.  Let $X'$ be homeomorphic to $X$ and parametrized as $T^2 \times [0,\infty)$, $T^2 \times \R$ or $S^1 \times \R^2$.  For each slope at infinity, there exist infinitely many tight contact structures on $X$ with that slope, distinct up to proper isotopy, which do not extend to a tight contact structure on $X'$.

\end{thm}

\noindent This result stands in contrast to Eliashberg's original examples, all of which are neighborhoods of a transverse curve in $S^3$ and have a different slope at infinity.  Using this embedding, it is easy to compactify his examples.  Theorem~\ref{thm: uncountableembeddings} shows that, in general, finding such a nice compactification is not straightforward.

Finally, just as high torus division number is a problem in the classification of toric annuli, contact structures with
infinite division number at infinity prove difficult to understand.  However, we are able to use the notion of {\em stable
disk equivalence} to partially understand this situation.  Precise statements of all of these results are in
Section~\ref{tighttoricends}.  In Section~\ref{opentori}, we use these results to reduce the classification of tight
contact structures on $S^1 \times \R^2$ and $T^2 \times \R$ to the classification of the corresponding toric ends.  

In the second half of the paper, we use the notion of the slope at infinity to prove a generalization of Eliashberg's
result in \cite{El3}: 

\begin{thm} \label{thm: uncountable}

Let $V$ be any open $3$-manifold which is the interior of a compact, irreducible, connected $3$-manifold $M$ such that
$\partial M$ is nonempty and contains at least one component of nonzero genus.  Then $V$ supports uncountably many tight
contact structures which are not contactomorphic, yet are isotopic.

\end{thm}

Eliashberg's proof involves computing the contact shape of the contact structures on $S^1 \times \R^2$, which in turn
relies on a previous computation of the symplectic shape of certain subsets of $T^n \times \R^n$ done in \cite{Si}.  We
bypass the technical difficulties of computing the symplectic shape by employing convex surface theory in the end of $V$.
The first step in the proof is to put a tight contact structure on the manifold $M$ with a certain dividing curve
configuration on the boundary.  To do this, we use the correspondence between taut sutured manifolds and tight contact
structures covered in \cite{HKM2}.  We then find nested sequences of surfaces which allow us to construct a contact manifold
$(V, \eta_s)$ for every $s \in (-2,-1)$.  We distinguish these contact structures up to proper isotopy by showing that they
have different slopes at infinity.  Since the mapping class group of an irreducible $3$-manifold with boundary is countable
(see \cite{McC}), uncountably many of the $\eta_s$ are not contactomorphic.  To simplify the presentation of the proof, we
first present the proof in the case when $\bdry M$ is connected in Section~\ref{incomp}.  We deal with
the case of disconnected boundary in Section~\ref{gencase}.    

In \cite{El1}, Eliashberg declares a contact structures on an open $3$-manifold $V$ to be {\em overtwisted at infinity} if
for every relatively compact $U \subset V$, each noncompact component of $V \setminus U$ is overtwisted.  If the contact
structure is tight outside of a compact set, then it is {\em tight at infinity}.  He then uses his classification for
overtwisted contact structures in \cite{El2} to show that any two contact structures that are overtwisted at infinity and
homotopic as plane fields are properly isotopic.  In contrast to this result, we have the following: 

\begin{thm} \label{thm: otuncountable} 

Let $V$ be any open $3$-manifold which is the interior of a compact, irreducible, connected $3$-manifold $M$ such that
$\partial M$ is nonempty and contains at least one component of nonzero genus.  Then $V$ supports uncountably many
overtwisted contact structures which are tight at infinity and which are not contactomorphic, yet are isotopic. 

\end{thm}


\section{Background and Conventions}


For general facts about $3$-manifolds, we refer the reader to \cite{He}.  For terminoloy and facts about contact geometry
and especially convex surface theory, we refer to \cite{Ho2} and \cite{Et}.  Given a convex surface $S$ in a contact
$3$-manifold, we denote the dividing set of $S$ by $\Gamma_S$.  The Legendrian Realization Principle
(see \cite{Ho2}) says that any nonisolating collection of arcs and closed curves on a convex surface can be made
Legendrian after an isotopy of the surface.  When we say ``LeRP'', we will mean ``apply the Legendrian Realization
Principle'' to a collection of curves.  We will use this as a verb and call this process ``LeRPing'' a collection of curves.

For the reader's convenience, we list some of the definitions and results in \cite{HKM2} which we will need later.  A {\em sutured manifold} $(M, \gamma)$ is a compact oriented $3$-manifold $M$ together with a set $\gamma \subset \partial
M$ of pairwise disjoint annuli $A(\gamma)$ and tori $T(\gamma)$.  $R(\gamma)$ denotes $\partial M \setminus
\intr(\gamma)$.  Each component of $R(\gamma)$ is oriented.  $R_+ (\gamma)$ is defined to be those components of
$R(\gamma)$ whose normal vectors point out of $M$ and $R_\gamma$ is defined to be $R(\gamma) \setminus R_+ (\gamma)$. 
Each component of $A(\gamma)$ contains a {\em suture} which is a homologically nontrivial, oriented simple closed curve. 
The set of sutures is denoted $s(\gamma)$.  The orientation on $R_+(\gamma)$, $R_-(\gamma)$ and $s(\gamma)$ are related as
follows.  If $\alpha \subset \partial M$  is an oriented arc with $\partial \alpha \subset R(\gamma)$ that intersects
$s(\gamma)$ transversely in a single point and if $s(\gamma) \cdot \alpha = 1$, then $\alpha$ must start in $R_+ (\gamma)$
and end in $R_- (\gamma)$ .

A {\em sutured manifold with annular sutures} is a sutured manifold $(M, \gamma)$ such that $\partial M$ is nonempty, every
component of $\gamma$ is an annulus, and each component of $\partial M$ contains a suture.  A sutured manifold $(M, \gamma)$
with annular sutures determines an {\em associated convex structure} $(M, \Gamma)$, where $\Gamma = s(\gamma)$.  For more on
this correspondence, see \cite{HKM2}.

A transversely oriented codimension-$1$ foliation $\mathcal{F}$ is {\em carried by} $(M, \gamma)$ if $\mathcal{F}$ is
transverse to $\gamma$ and tangent to $R(\gamma)$ with the normal direction pointing outward along $R_+(\gamma)$ and inward
along $R_-(\gamma)$, and $\mathcal{F}|_{\gamma}$ has no Reeb components.  $\mathcal{F}$ is {\em taut} if each leaf
intersects some closed curve or properly embedded arc connecting $R_-(\gamma)$ to $R_+(\gamma)$ transversely.

Let $S$ be a compact oriented surface with components $S_1, \dots, S_n$.  Let $\chi(S_i)$ be the Euler characteristic of
$S_i$.  The {\em Thurston norm of} $S$ is defined to be $$x(S) =  \sum_{\chi(S_i)<0} |\chi(S_i)|.$$

A sutured manifold $(M, \gamma)$ is {\em taut} if 

\begin{enumerate}
\item $M$ is irreducible.
\item $R(\gamma)$ is Thurston norm minimizing in $H_2(M, \gamma)$; that is, if $S$ is any other properly embedded surface
with $[S] = [R(\gamma)]$, then $x(R(\gamma)) \leq x(S)$.
\item $R(\gamma)$ is incompressible in $M$.
\end{enumerate}

The following is due to Gabai \cite{Ga} and Thurston \cite{Th}.

\begin{thm} 

A sutured manifold $(M, \gamma)$ is taut if and only if it carries a transversely oriented, taut, codimension-$1$ foliation
$\mathcal{F}$.

\end{thm}

We require the following result due to Honda, Kazez, and Mati\'c \cite{HKM2}.

\begin{thm}
Let $(M, \gamma)$ be an irreducible sutured manifold with annular sutures, and let $(M, \Gamma)$ be the associated convex
structure.  The following are equivalent.

\begin{enumerate}

\item $(M, \gamma)$ is taut.
\item $(M, \gamma)$ carries a taut foliation.
\item $(M, \Gamma)$ carries a universally tight contact strucuture.
\item $(M, \Gamma)$ carries a tight contact structure.

\end{enumerate}

\end{thm}


\section{The end of an open contact manifold and some invariants}\label{invariant} 


Let $(V, \xi)$ be any open contact $3$-manifold which is the interior of a compact $3$-manifold $M$ such that $\bdry M$ is
nonempty and contains at least one component of nonzero genus.  Fix an embedding of $V \hookrightarrow int(M)$ so that we
can think of $V$ as $M \setminus \partial M$.  Choose a boundary component $S \subset \partial M$ and let $\Sigma
\subset M \setminus \bdry M$ be an embedded surface isotopic to $S$ in $M$.  Note that $S$ and $\Sigma$ bound a contact manifold $(\Sigma
\times (0,1), \xi)$.  We call such a manifold, along with the embedding into $V$, a {\em contact end} corresponding to $S$
and $\xi$.  Let $Ends(V, \xi; S)$ be the collection of contact ends corresponding to $S$ and $\xi$.  

Let $S \subset \partial M$ be a component of nonzero genus and let $\lambda \subset S$ be a separating, simple closed curve
which bounds a punctured torus $T$ in $S$.  Fix a basis $B$ of the first homology of $T$.  Let $\Sigma \subset
V$ be a convex surface which is isotopic to $S$ in $M$ and contains a simple closed curve $\gamma$ with the following properties:

\be  

\item $\gamma$ is isotopic to $\lambda$ on $\Sg$, where we have identified $\Sg$ and $S$ by an isotopy in $M$.
\item $\gamma$ intersects $\Gamma_{\Sigma}$ transversely in exactly two points.
\item $\gamma$ has minimal geometric intersection number with $\Gamma_{\Sigma}$.

\ee

Call any such surface {\em well-behaved} with respect to $S$ and $\lambda$.  Note that there exists a simple closed curve
$\mu \subset \Gamma_{\Sigma}$ which is contained entirely in $T$.  Let the {\em slope} of $\Sigma$, written $\slope(\Sg)$,
be the slope of $\mu$ measured with respect to the basis $B$ of the first homology of $T$.  When $S$ is a torus, we omit
all reference to the curve $\lambda$ as it is unnecessary for our definition.

Let $E \in Ends(V, \xi; S)$.  Let $\CC(E)$ be the set of all well-behaved convex surfaces in the contact end $E$.  If
$\CC(E) \neq \emptyset$, then define the {\em slope} of $E$, to be $$slope(E) \; = \; \sup_{\Sg \in \CC(E)}(\slope(\Sg)).$$ 
Here we allow $\text{sup}$ to take values in $\R \cup \{\infty\}$.  Note that $Ends(V, \xi; S)$ is a directed set, directed
by reverse inclusion and that the function $\slope \co Ends(V, \xi; S) \rightarrow \R \cup \infty$ is a net.  If $\CC(E)$
is nonempty for a cofinal sequence of contact ends and this net is convergent, then we call the limit the {\em slope at
infinity} of $(V, \xi; S, \lambda, B)$ or the {\em slope at infinity} of $(V, \xi)$ if $S$, $\lambda$, and $B$ are
understood from the context.  If the slope at infinity exists, then we say that this slope is {\em attained} if for each $E
\in Ends(V, \xi; S)$ there exists a $\Sg \in \CC(E)$ with that slope.  Note that any slope that is attained must
necessarily be rational.

Let $\Sg \in \CC(E)$.  Define the {\em division number} of $\Sg$, written $\dv(\Sg)$ to be half the number of dividing
curves and arcs on $T$.  When $\Sg$ is a torus, this is the usual torus division number.  If $\CC(E) \neq \emptyset$, then
let $$\dv(E) = \min_{\Sg \in \CC(E)}(\dv(\Sg)).$$  Note that $\dv \co Ends(V, \xi, S) \ra \N \cup \{ \infty \}$ is a net,
where we endow $\N \cup \{ \infty \}$ with the discrete topology.  If $\CC(E)$ is nonempty for a cofinal sequence of
contact ends, then we call the limit the {\em division number at infinity} of $(V, \xi; S, \lambda, B)$ or the {\em
division number at infinity} of $(V, \xi)$ if $S$, $\lambda$, and $B$ are understood from the context.  Note that the slope
at infinity and the division number at infinity are proper isotopy invariants.


\section{Classification Theorems for Tight Toric Ends} \label{tighttoricends}


In this section, we study tight contact structures on toric ends.  We say that a toric end is {\em minimally twisting} if
it contains only minimally twisting toric annuli.  We first show that it is possible to refer to the slope at infinity and
the division number at infinity for toric ends.

\begin{prop}

Let $T^2 \times [0,\infty)$ be a tight toric end.  Then the division number at infinity and the slope at infinity are
defined.

\end{prop}

\begin{proof} 

First note that $\CC(E)$ is nonempty for any end $E$ since the condition for being well-behaved is vacuously
true for tori. Also, note that the division number at infinity exists by definition.

If there exists a nested sequence of ends $E_i$ such that $\slope(E_i) = \infty$, then the slope at infinity is $\infty$. 
Otherwise, there exists an end $E = T^2 \times [0, \infty)$ such that for no end $F \subset E$ is $\slope(F) = \infty$. 
This means that $E$ is minimally twisting.  Without loss of generality, assume $T_i = T^2 \times i$ is convex with slope $s_i$. 
Note that the $s_i$ form a clockwise sequence on the Farey graph and are contained in a half-open arc which does not
contain $\infty$.  Since $slope(F) \leq s_i$ for any end $F \subset T^2 \times [i, \infty)$, our net is convergent, so the
slope at infinity is defined. \end{proof}


\subsection{Tight, minimally twisting toric ends with irrational slope at infinity} 


In this section, we study tight, minimally twisting toric ends $(T^2 \times [0,\infty), \xi)$ with {\em irrational} slope
$r$ at infinity and with convex boundary satisfying $\dv(T^2 \times 0) = 1$ and $\slope(T^2 \times 0) = -1$.  Unless
otherwise specified, all toric ends will be of this type.  

We first show how to associate to any such toric end a function $f_{\xi} \co \N \ra \N \cup \{0\}$.  There exists a sequence of
rational numbers $q_i$ on the Farey graph which satisfies the following:  

\be

\item $q_1 = -1$ and the $q_i$ proceed in a clockwise fashion on the Farey graph. 
\item $q_i$ is connected to $q_{i+1}$ by an arc of the graph. 
\item The $q_i$ converge to $r$.  
\item The sequence is minimal in the sense that $q_i$ and $q_j$ are not joined by an arc of the graph unless $j$ is
adjacent to $i$.

\ee

We can form this sequence inductively by taking $q_2$ to be the rational number which is closest to $r$ on the clockwise
arc of the Farey graph $[-1,r]$ between $-1$ and $r$ and has an edge of the graph from $-1$ to $q_2$.  Similarly, construct
the remaining $q_i$.  Any such sequence can be grouped into {\em continued fraction blocks}.  We say that $q_i, \dots, q_j$
form a {\em continued fraction block} if there is an element of $SL_2(\Z)$ taking the sequence to $-1,\dots,-m$.  We call
$m$ the {\em length} of the continued fraction block.  We say that this block is {\em maximal} if it cannot be extended to a
longer continued fraction block in the sequence $q_i$.  Since $r$ is irrational, maximal continued fraction blocks exist. 
Denote these blocks by $B_i$.  To apply this to our situation, we need the following.

\begin{prop} 

There exists a nested sequence of convex tori $T_i$ with $\dv(T_i) = 1$ such that $\slope(T_i) = q_i$. 
Moreover, any such sequence must leave every compact set. 

\end{prop}

\begin{proof}

By the definition of slope at infinity, for any $\epsilon$, there is an end $E$ such that $\slope(E)$ is within $\epsilon$
of $r$.  This means that there is a convex torus $T$ in $E$ with slope lying within $2\epsilon$ of $r$.  Note that since
our toric end is minimally twisting and has slope $r$ at infinity, $\slope(T) \in [-1, r)$.  We attach bypasses to $T$
so that $\dv(T)=1$.  The toric annulus bounded by $T^2 \times 0$ and $T$ contains the tori $T_i$ with $q_i$ lying
couterclockwise to $\slope(T)$.  Fix these first $T_i$.  Choose another torus $T'$ outside of the toric annulus with
slope even closer to $r$.  Again, adjust the division number of $T'$ so that it is $1$ and factor the toric annulus bounded
by $T$ and $T'$ to find another finite number of our $T_i$.  Proceeding in this fashion, we see we have the desired
sequence of $T_i$.  Any such sequence must leave every compact set by the definition of the slope at infinity.  For,
if not, then we could find a torus $T$ in any end with $\slope(T) > r$, which would show that the slope at infinity is not
$r$. \end{proof}

This factors the toric end according to our sequence of rationals.  We say that a consecutive sequence of $T_i$ form a
continued fraction block if the corresponding sequence of rationals do.  Each maximal continued fraction block $B_i$
determines a maximal continued fraction block of tori which we also call $B_i$.  We think of $B_i$ as a toric annulus.

To each continued fraction block, we let $n_j$ be the number of positive basic slices in the factorization of $B_i$ by
$T_j$.  Define $f_\xi \co \N \ra \N \cup \{0\}$ by $f_\xi(j) = n_j$.  To show that the function $f_\xi$ is independent of
the factorization by $T_i$, suppose $T_i'$  is another factorization with the same properties as $T_i$.  Let $B_j'$ denote
the corresponding continued fraction blocks.  Fix $j$.  There exists $n$ large such that the toric annulus $A$ bounded by
$T_n$ and $T_1$ contains the continued fraction blocks $B_j$ and $B_j'$.  Extend the partial factorization of $A$ by
$B_j'$.  Recall that one can compute the relative Euler class via such a factorization and that it depends on the number of
positive basic slices in each continued fraction block \cite{Ho2}.  Therefore, $B_j$ and $B_j'$ must have the same number of positive
basic slices.

Given an irrational number $r$, let $\CF(r)$ denote the collection of functions $f \co \N \ra \N \cup \{0\}$ such that
$f(i)$ does not exceed one less than the length of $B_i$.  We can now state a complete classification of the toric ends
under consideration.

\begin{thm} \label{thm: mintwistirr}

Let $(T^2 \times [0, \infty), \xi)$ be a tight, minimally twisting toric end with convex boundary satisfying $\dv(T^2
\times 0) = 1$ and $\slope(T^2 \times 0) = -1$.  Suppose that the slope at infinity is irrational.  To each such tight
contact structure, we can assign a function $f_{\xi}\co \N \ra \N \cup \{0\}$ which is a complete proper isotopy (relative
to the boundary) invariant.  Moreover, given any $f \in \CF(r)$, there exists a toric end $(T^2 \times [0,\infty), \xi)$
such that $f_\xi = f$.

\end{thm}

\bp

If $f_\xi = f_{\xi'}$, then we can shuffle bypasses within any given continued fraction block so that all positive basic
slices occur at the beginning of the block.  Since the number of positive basic slices in any continued fraction block is
the same, it is clear that they are properly isotopic.

It is a straightforward application of the gluing theorem for basic slices in \cite{Ho2} to show that we can construct a
toric annulus corresponding to the desired continued fraction blocks.  The fact that they stay tight under gluing follows
from the fact that overtwisted disks are compact.\ep

\begin{cor} \label{cor: irrembeddings}

Let $(T^2 \times [0,1),\xi)$ be a tight, minimally twisting toric end with irrational slope $r$ at infinity.  Suppose $f_\xi(i)$ is not maximal or minimal for an infinite number of numbers $i$.   Then there does not exist any tight, toric end $(T^2 \times [0,\infty), \eta)$ such that $\xi|_{T^2 \times [0,1)} = \eta|_{T^2 \times [0,1)}$.  \end{cor}

\bp

Assume that there were an inclusion $\phi \co (T^2 \times [0,1),\xi) \ra (T^2 \times [0,\infty),\eta)$.  Perturb $T^2 \times \{2\}$ to be convex of slope $b$.  Choose a convex torus $\phi(T')$ of slope $a$.  As before, we have a minimal, clockwise sequence of rationals $q_j$ for $1 \leq j \leq n$ on the Farey graph such that $q_1 = a$, $q_n = b$, and $q_i$ is joined to $q_{i+1}$ by an arc of the graph.  Let $q_m$ be the rational closest to $q_1$ such that $r$ lies clockwise to $q_1$ and counterclockwise to $q_m$.  By our assumption on $f_\xi$, there exists a continued fraction block of tori $T_{j_1}, \dots , T_{j_k} \subset (T^2 \times [0,1), \xi)$ which contains both positive and negative basic slices.  Moreover, we can assume that the corresponding sequence of rationals lies clockwise to $q_{m-1}$ and counterclockwise to $q_m$.  Perturb tori $T_{in}$ and $T_{out}$ in $(T^2 \times [0,\infty),\eta)$ to be convex of slopes $q_{m-1}$ and $q_m$, respectively, such that the basic slice bounded by $T_{in}$ and $T_{out}$ contains $\phi(T_{j_1}), \dots , \phi(T_{j_k})$.  This is a contradiction, since a basic slice cannot be formed by gluing basic slices of opposite signs unless the contact structure $\eta$ is overtwisted \cite{Ho2}.\ep


\subsection{Tight, minimally twisting toric ends with rational slope at infinity} 


We now consider tight, minimally twisting toric ends $(T^2 \times [0,\infty), \xi)$ with {\em rational} slope $r$ at
infinity and with convex boundary satisfying $\dv(T^2 \times 0) = 1$ and $\slope(T^2 \times 0) = -1$.  Unless otherwise
specified, all toric ends will be of this type.  We first deal with the situation when the slope at infinity is {\em not}
attained.  

We show how to every toric end under consideration we can assign a function $f_{\xi} \co \{1, \dots, n(r)\} \times \{1, -1
\} \ra \N \cup\{0, \infty\}$.  We proceed in a fashion similar to the irrational case.  Given $r$ rational, there exists a
sequence of rationals $q_i$ satisfying the following:

\be

\item $q_1 = -1$ and the $q_i$ proceed in a clockwise fashion on the Farey graph.
\item $q_i$ is connected to $q_{i+1}$ by an arc of the tesselation.
\item The $q_i$ converge to $r$, but $q_i \neq r$ for any $i$. 
\item The sequence is minimal in the sense that $q_i$ and $q_j$ are not joined by an arc of the tesselation unless $j$ is
adjacent to $i$.

\ee

We construct such a sequence inductively just as in the irrational case, except we never allow the rationals $q_i$ to reach
$r$.  Note that such a sequence breaks up naturally into $n-1$ finite continued fraction blocks $B_i$ and one infinite
continued fraction block $B_{n}$ (\ie, $B_{n}$ can be taken to the negative integers after action by $SL_2(\Z)$).  Note
that $n$ is completely determined by $r$.  Just as in the irrational case, there exist nested covex tori $T_i$ with
$\dv(T_i)=1$ and $\slope(T_i) = q_i$.  We can argue as in the irrational case to show that these tori must leave every
compact set of the toric end.  We will also refer to the collection of tori $T_i$ corresponding to $B_i$ by the same name.

We will now construct $f_\xi$.  Let $f_\xi(i, \pm 1)$ be the number of positive (negative) basic slices in the continued
fraction block $B_i$.  Of course, for a finite continued fraction block, $f_\xi(i, 1)$ determines $f_\xi(i,-1)$.  However,
this is clearly not the case for $B_n$.

As in the irrational case, let $\CF(r)$ be the collection of functions $f \co \{1, \dots, n(r)\} \times \{1, -1 \} \ra \N
\cup\{0, \infty\}$ such that $f_\xi(i, 1) + f_\xi(i,-1) = |B_i| - 1$ for $i \leq n-1$, where $|B_i|$ is the length of
$B_i$, and at least one of $f_\xi(n(r),\pm 1)$ is infinite.

\begin{thm} \label{thm: mintwistrat}

Let $(T^2 \times [0,\infty), \xi)$ be a tight, minimally twisting toric end with convex boundary satisfying $\dv(T^2 \times
0) = 1$ and $\slope(T^2 \times 0) = -1$.  Suppose that the slope at inifinity is rational and is not attained.  To each
such tight contact structure, we can assign a function $f_{\xi}: \{1, \dots, n(r)\} \times \{1, -1 \} \ra \N \cup\{0,
\infty\}$ which is a complete proper isotopy (relative to the boundary) invariant.   Moreover, for any $f \in \CF(r)$,
there exists a tight, minimally twisting toric end $(T^2 \times [0,\infty), \xi)$ with slope $r$ at infinity which is not
realized such that $f = f_\xi$.

\end{thm}

\bp

Suppose $f_\xi = f_\xi'$.  As in the irrational case, we can adjust our factorization of the finite continued fraction
blocks so that all of the positive basic slices occur first in each continued fraction block.  Therefore, we can isotope the
two contact structures so that they agree on the first $n-1$ continued fraction blocks.  

We now consider the infinite basic slice.  Without loss of generality, we may assume that the infinite basic slices for
$\xi$ and $\xi'$ are toric ends $(T^2 \times [0,\infty), \xi)$ and $(T^2 \times [0,\infty), \xi')$ with $\slope(T^2 \times
\{0\})$, $\dv(T^2 \times \{0\})$, and infinite slope at infinity that is not realized.  The corresponding factorization is
then given by nested tori $T_i$ and $T_i'$ such that $\slope(T_i) = \slope(T_i') = -i$ and $\dv(T_i) = 1$.  We now construct
model toric ends $\xi_n^{\pm}$ and $\xi_{alt}$ and show that any infinite basic slice is properly isotopic to one of the
models.  Let $B_i^{\pm}$ be the positive (negative) basic slice with $\slope(T^2 \times 0) = -i$ and $\slope(T^2 \times 1)
= -i-1$.  Let $\xi_n^{\pm}$ be the toric end constructed as $B_1^{\pm} \cup \dots \cup B_n^{\pm} \cup B_{n+1}^{\mp} \cup
\cdots$.  Let $\xi_{alt}$ be $B_1^+ \cup B_2^- \cup B_3^+ \cup \cdots$.  First consider the case when $f_\xi(n,1) = m$. 
There exists $N$ large so that the toric annulus bounded by $T_1$ and $T_N$ contains {\em at least} $m$ positive basic
slices and $m$ negative basic slices.  By shuffling bypasses in this toric annulus, we can rechoose our factorization so
that all positive bypass layers occur first in our factorization.  This toric end is clearly properly isotopic to
$\xi_m^+$.  We handle the case when $f_\xi(n,-1) = m$ similarly.  Now, suppose that $f_\xi(n, \pm 1) = \infty$.  Fix some
number $k$.  Choose $N_1$ large enough that the toric annulus bounded by $T_1$ and $T_{N_1}$ contains at least $k$ positive
and $k$ negative basic slices.  By shuffling bypasses in this toric annulus, we can arrange for the first $2k$ basic slices
in the factorization to be alternating.  There exists an isotopy $\phi_t^1$ such that $\phi_0^1$ is the identity and
${\phi_{1}^1}_*(\xi)$ agrees with $\xi_{alt}$ in the first $2k$ basic slices.  Call the pushed forward contact structure by
the same name.  There exists $N_2$ large such that $T_{2k}$ and $T_{N_2}$ bound a toric annulus with $k$ positive and $k$
negative basic slices.  Leaving the first $2k$ tori in our factorization fixed, we can shuffle bypases in the toric annulus
bounded by $T_{2k}$ and $T_{N_2}$ so that signs are alternating.  Choose an isotopy $\phi_t^2$ as before such that
$\phi_t^2$ is the identity on the toric annulus bounded by $T_1$ and $T_{2k}$ and takes the second $2k$ basic slices of
$\xi$ onto those of $\xi_{alt}$.  Continuing in this fashion, we can construct $\phi_t^n$ which is supported on $K_n$
compact such that $K_i \subset K_{i+1}$ and $T^2 \times [0,\infty) = \cup K_i$.  Hence we have an isotopy taking $\xi$ to
$\xi_{alt}$.  The existence result follows immediately from Honda's gluing results for toric annuli \cite{Ho2}.\ep

\begin{cor} \label{cor: ratembeddings}

Let $(T^2 \times [0,1),\xi)$ be a tight, minimally twisting toric end that does not attain a rational slope $r$ at infinity.  Suppose
$f_\xi(n(r) \times \{1\})$ and $f_\xi(n(r) \times \{-1\})$ are nonzero.  Then there does not exist any tight, toric end $(T^2 \times [0,\infty), \eta)$ such that $\xi|_{T^2 \times [0,1)} = \eta|_{T^2 \times [0,1)}$. \end{cor}

\bp

Assume that there were such an inclusion  $\phi \co (T^2 \times [0,1),\xi) \ra (T^2 \times [0,\infty),\eta)$.  Let $T_i$ be the first torus in the factorization of the infinite continued fraction block of $(T^2 \times [0,1),\xi)$.  By definition, there exists another torus $T_j$ with $j > i$ such that $T_i$ and $T_j$ bound basic slices of both signs.  By the definition of the slope at infinity and the precompactness condition, there exists a convex torus $T$ outside of the toric annulus bounded by $\phi(T_i)$ and $\phi(T_j)$ which has slope $r$.  Note that $\phi(T_i)$ and $T$ bound a continued fraction block which is formed by gluing basic slices of opposite signs.  This implies that $(T^2 \times [0,\infty),\eta)$ is overtwisted \cite{Ho2}.\ep

Corollary~\ref{cor: irrembeddings} and Corollary~\ref{cor: ratembeddings} will be essential to proving Theorem~\ref{thm: uncountableembeddings}.  We now consider tight, minimally twisting toric ends that realize the slope at infinity and have finite division number at
infinity.

\begin{thm} \label{thm: mintwistratach}

Tight, minimally twisting toric ends with finite division number $d$ at infinity that realize the slope $r$ at infinity are
in one-to-one correspondence with tight, minimally twisting contact structures on $T^2 \times [0,1]$ with $T^2 \times i$
convex, $\slope(T^2 \times 0) = -1$, $\slope(T^2 \times 1) = r$, $\dv(T^2 \times 0) = 1$, and $\dv(T^2 \times 1) = d$ up to
isotopy relative to $T^2 \times 0$. 

\end{thm} 

\bp

Let $(T^2 \times [0,\infty), \xi)$ be such a toric end.  By the definition of division number at infinity and slope at
infinity, there exists a convex torus $T$ with the following properties:

\be

\item $\dv(T) = d$
\item $\slope(T) = r$
\item Any other convex torus $T'$ lying in the noncompact component of $T^2 \times [0,\infty) \setminus T$ satisfies
$\dv(T') \geq d$.  

\ee

Any such torus will necessarily have slope $r$.  Let $A$ be the toric annulus bounded by $T^2 \times 0$ and $T$.  We know
that any other torus $T'$ with the same properties as $T$ and bounds a toric annulus $A'$ is isotopic to $A$.  By the
definition of $T$ and $T'$ there exists a torus $T''$ outside of $A$ and $A'$ that has the same properties as $T$.  Since
$\xi$ is minimally twisting, $T'$ and $T''$ bound a vertically invariant toric annulus.  Similarly, $T$ and $T''$ bound a
vertically invariant toric annulus.  We can use these toric annuli to isotope $A$ and $A'$ to the same toric annulus in our
toric end.  This yields the desired correspondence.  Given a tight, minimally twisting contact structures on $T^2 \times
[0,1]$ with $T^2 \times i$ convex, $\slope(T^2 \times 0) = -1$, $\slope(T^2 \times 1) = r$, $\dv(T^2 \times 0) = 1$, and
$\dv(T^2 \times 1) = d$, we obtain a toric end by removing $T^2 \times 1$.\ep

We say that two convex annuli $A_i = S^1 \times [0,1]$ with Legendrian boundary, $\tb(S^1 \times 0) = -1$ and $\tb(S^1
\times 1) = -m$ are {\em stabily disk equivalent} if there exist disk equivalent convex annuli $A_i' = S^1 \times [0,2]$
such that $\tb(S^1 \times 1) = -1$, $\tb(S^1 \times 2) = -n < -m$, and $A_i = S^1 \times [0,1] \subset A_i'$.

\begin{thm} \label{thm: infdivnum}

Let $(T^2 \times [0,\infty), \xi)$ be a tight, minimally twisting toric end with $\slope(T^2 \times 0) = \infty$, slope
$\infty$ at infinity, and division number $\infty$ at infinity.  Then we can associate to $\xi$ a collection of nested
families of convex annuli $A_i = S^1 \times [0, i]$ with Legendrian boundary such that $\tb(S^1 \times 0) = -1$, $\tb(S^1
\times i+1) = \tb(S^1 \times i)+1$ such that any two annuli $A_i$ and $A_i'$ in different families are stabily disk
equivalent.

\end{thm}

\bp

To construct such annuli, simply choose a factorization of the toric end by tori $T_i$ such that $T_1 = T^2 \times 0$,
$\slope(T_i) = \infty$, $\dv(T_{i+1}) = \dv(T_i)+1$ and the $T_i$ leave every compact set.  Let $A_1$ be the convex annulus
with boundary on $T_1$ and $T_2$.  Choose $A_1'$ a horizontal convex annulus between $T_2$ and $T_3$ which shares a
boundary component with $A_1$.  Let $A_2 = A_1 \cup A_1'$.  Continuing in this fashion, we construct a sequence of nested
annuli $A_i$.  Now, choose any other factorization by tori $T_i'$ satisfying the same properties as the $T_i$ and let
$A_i'$ be the corresponding sequence of convex annuli.  We will show that $A_i$ is stabily disk equivalent to $A_i'$. 
Choose $N$ large so that the toric annulus bounded by $T_1$ and $T_N$ contains $A_i$ and $A_i'$.  Let $A$ be a convex
annulus between the $S^1 \times i \subset A_i'$ and a horizontal Legendrian curve on $T_N$.  Let $A' = A_i' \cup A$. 
Honda's result in \cite{Ho2} implies that $A$ and $A'$ are disk equivalent.\ep

\begin{cor} \label{cor: infdivnumemb}

Any tight, minimally twisting toric end $(T^2 \times [0,\infty), \xi)$ with $\slope(T^2 \times 0) = \infty$, slope $\infty$
at infinity, and division number $\infty$ at infinity embeds in a vertically invariant neighborhood of $T^2 \times 0$.

\end{cor}

\bp

Honda's model \cite{Ho2} for increasing the torus division number can be applied inductively on a vertically invariant neighborhood of
$T^2 \times 0$ to create the desired sequence of nested tori $T_i$ and corresponding annuli $A_i$.  The contact structure on
the toric annulus bounded by $T_1$ and $T_i$ is uniquely determined by $A_i$ \cite{Ho2}.\ep

\noindent We are lead to the following question:

\begin{question}

What are necessary and sufficient conditions for two toric ends with infinite division number at infinity to be properly
isotopic?

\end{question}


\subsection{Nonminimally twisting, tight toric ends}


In this section, we deal with tight toric ends $(T^2 \times [0, \infty), \xi)$ with $\slope(T^2 \times 0) = 0$, $\dv(T^2
\times 0) = 1$, and are not minimally twisting.  We first recall Honda's classification for nonminimally twisting tight
contact structures on $T^2 \times [0,1]$ in \cite{Ho2}.  He constructs a family $\xi_n^{\pm}$ of tight, rotative contact
structures on $T^2 \times [0,1]$ with $\slope(T^2 \times i) = 0$ and $\dv(T^2 \times i) = 1$ and shows that this is a
complete and nonoverlapping list of contact structures satisfying these conditions.  We define the {\em rotativity} of a
tight toric end $\xi$ with $\slope(T^2 \times 0) = 0$ and $\dv(T^2 \times 0) =1$ to be the maximum $n$ such that there is
an embedding $e \co (T^2 \times [0,1], \xi_n^{\pm}) \hra (T^2 \times [0,\infty), \xi)$ with $e(T^2 \times 0) = T^2 \times
0$.  If no maximum exists, then we say that $\xi$ has {\em infinite rotativity}.  If $n$ is the rotativity of $\xi$, then
$\xi_n^{+}$ and $\xi_n^{-}$ cannot both be embedded in $\xi$.  For, the images of any two such embeddings would provide two
factorizations for a common toric annulus.  But, such factorizations are unique \cite{Ho2}.  Hence, we can refer to the
{\em sign of rotativity} as well.  We construct two more nonminimally twisting toric ends $\xi_{\infty}^{\pm}$.  Set $(T^2
\times [0, \infty), \xi_{\infty}^{\pm}) = \cup_{i=1}^{\infty}(T^2 \times [0,1], \xi_2^{\pm})$.

\begin{thm} \label{thm: nonminfac}

Let $(T^2 \times [0, \infty), \xi)$ be a tight toric end which is not minimally twisting.  

\be

\item Assume that $\xi$ has finite rotativity $n$ and that the sign of rotativity is $+$.  Let $e, e' \co (T^2 \times [0,1],
\xi_n^{\pm}) \hra (T^2 \times [0,\infty), \xi)$ be any two embeddings with $e(T^2 \times 0) = e'(T^2 \times 0) = T^2 \times
0$.  Let $T = e(T^2 \times 1)$ and $T^\prime = e^\prime(T^2 \times 1)$.  Then the contact structures on the toric ends
bounded by $T$ and $T^\prime$ are identical.  Moreover, $\xi$ is universally tight.   

\item Assume that $\xi$ has infinite rotativity.  Then $\xi$ is properly isotopic relative to the boundary to either
$\xi_{\infty}^+$ or $\xi_{\infty}^-$, so the sign of rotativity is defined in the infinite case as well.  Moreover, the 
$\xi_{\infty}^{\pm}$ are universally tight.

\ee

\end{thm}

\bp

First, consider the case of finite rotativity.  Both $T$ and $T'$ are contained in a toric annulus bounded by $T^2 \times 0$
and some torus $T''$ with nonzero slope and division number $1$.  The contact structures on the toric annulus $A$
($A^{\prime}$) bounded by $T$ and $T''$ ($T^\prime$ and $T^{\prime\prime}$) are uniquely determined by Honda's work
\cite{Ho2}.  Therefore, we can isotope the two factorizations so that they coincide.  The fact that $\xi$ is universally
tight follows immediately, since nonminimally twisting toric annuli are universally tight. 

Now, assume $\xi$ has infinite rotativity.  First, note that we cannot have two embeddings $e_n^{\pm} \co (T^2 \times [0,1],
\xi_n^{\pm}) \hra (T^2 \times [0,\infty), \xi)$ with $e_n^{\pm}(T^2 \times 0) = T^2 \times 0$ by the uniqueness of
factorizations of toric annuli.  Since $\xi$ has infinite rotativity, there exists a sequence of, say, positive embeddings
$e_n \co (T^2 \times [0,n], \xi_n^{+}) \hra (T^2 \times [0,\infty), \xi)$ with $e_n(T^2 \times 0) = T^2 \times 0$.  Moreover,
we can take this sequence of embeddings to be nested in the sense that $e_n = e_{n+1}$ on $[0,n]$.  This follows
immediately by factoring a toric annulus containing the images of $e_n$ and $e_{n+1}$.  Note that any sequence of such
embeddings must necessarily leave any compact set.  We can use this sequence of embeddings to construct a proper isotopy of
$\xi$ with $\xi_{\infty}^+$ as in the proof of Theorem~\ref{thm: mintwistrat}.  Again, the fact that $\xi_{\infty}^{\pm}$ are universally tight follows from the fact that
nonminimally twisting toric annuli are universally tight.\ep

\begin{cor} \label{cor: nonminembed}

Let $(T^2 \times [0, \infty), \xi)$ be a tight toric end that is not minimally twisting.  Then $(T^2 \times [0, \infty),
\xi)$ embeds into a toric annulus with convex boundary.  Moreover, if $\xi$ has finite division number at infinity, then
the image of the embedding can be chosen to have compact image.

\end{cor}


\section{Classifying Tight Contact Structures on $S^1 \times \R^2$ and $T^2 \times \R$} \label{opentori}


We now show that in many cases the classification of tight contact structures on $S^1 \times \R^2$ and $T^2 \times \R$
reduces to the classification of toric ends.


\subsection{Factoring tight contact structures on $S^1 \times \R^2$}


Let $(S^1 \times \R^2, \xi)$ be a tight contact structure and let $r$ be the slope at infinity.  Consider the collection of
points on the Farey graph of the form $1/n$ where $n \in \Z$.  Let $s(r) = 1/n$ be the point closest to $r$ (when
traversing the Farey graph counterclockwise) that is realized as the slope of a convex torus $T$ topologically isotopic to
$S^1 \times S^1$.  We can then factor $(S^1 \times \R^2, \xi)$ into $(S^1 \times D^2, \xi)$ and $(T^2 \times [0,\infty),
\xi)$.  To see that this factorization is unique, consider any other torus $T'$ satisfying the same conditions as $T$. 
Both $T$ and $T'$ lie in a common solid torus $S$ with convex boundary.  Note that the toric annuli bounded by $\bdry S$
and $T$ and by $\bdry S$ and $T^{\prime}$ are identical by the uniqueness of such factorizations on solid tori.  This
proves the following:

\begin{thm} \label{thm: solidtori}

Tight contact structures on $(S^1 \times \R^2, \xi)$ with nonzero slope at infinity are in one-to-one correspondence with
isotopy classes relative to the boundary of tight, minimally twisting toric ends $(T^2 \times [0,\infty), \eta)$.  Tight
contact structures on $(S^1 \times \R^2, \xi)$ with slope zero at infinity are in one-to-one correspondence with isotopy
classes relative to the boundary of tight, minimally twisting toric ends $(T^2 \times [0,\infty), \eta)$ which do not
attain the slope at infinity. 

\end{thm}  


\subsection{Factoring tight contact structures on $T^2 \times \R$}


In this section, we deal with tight contact structures on $T^2 \times \R$.  Any convex, incompressible torus $T \subset T^2
\times \R$ produces a factorization of $T^2 \times \R$ into $T^2 \times (-\infty, 0]$ and $T^2 \times [0,\infty)$.  We
identify $T^2 \times (-\infty, 0]$ with $T^2 \times [0,\infty)$ via reflection about the origin in $\R$ to obtain a
negative contact structure on $T^2 \times [0,\infty)$.  We change this to a positive contact structure by reflecting across
the $(1,0)$ curve in $T^2$.  Let $(T^2 \times [0,\infty), \xi_{\pm})$ and $(T^2 \times [0,\infty), \xi^{\prime}_{\pm})$ be
two factorizations corresponding to two different convex tori $T$ and $T^{\prime}$ with division number $1$ and slope $s$. 
We see that by keeping track of the $I$-twisting of a toric annulus in $T^2 \times \R$ containing $T$ and $T^\prime$, we
can obtain $(T^2 \times [0,\infty), \xi_{\pm})$ from $(T^2 \times [0,\infty), \xi^{\prime}_{\pm})$ as follows: Remove a
(possibly) rotative $T^2\times [0,1]$ with $\dv(T^2 \times i) = 1$ and $\slope(T^2 \times i) = s$ from the boundary of
$(T^2 \times [0, \infty), \xi_+)$ (or $(T^2 \times [0, \infty), \xi_-)$).  Apply a suitable diffeomorphism to $T^2 \times
[0,1]$.  Then, glue $T^2 \times [0,1]$ to the boundary of $(T^2 \times [0, \infty), \xi_-)$ (or $(T^2 \times [0, \infty),
\xi_+)$).  We call this procedure {\em shifting the rotativity} between $(T^2 \times [0,\infty), \xi_+)$ and $(T^2 \times
[0, \infty), \xi_-)$.

\begin{thm} \label{thm: opentoricannuli}

Let $(T^2 \times \R, \xi)$ be a tight contact manifold which contains a convex, incompressible torus $T$ with $\dv(T) = 1$
and $slope(T) = s$.  Then the factorization of $(T^2 \times \R, \xi)$ into toric ends $(T^2 \times [0,\infty), \xi_{\pm})$
is unique up to shifting the rotativity between the two toric ends.

\end{thm}

Theorem~\ref{thm: opentoricannuli} shows that the classification of contact structures on $T^2 \times \R$ reduces to the
study of toric ends if there is a convex, incompressible torus $T$ with $\dv(T) = 1$.  If $(T^2 \times \R, \xi)$ contains
no such torus, then the situation is much more subtle.

\begin{question}

If $(T^2 \times \R, \xi)$ contains no convex, incompressible torus with division number $1$, then what is the relationship
between two factorizations by convex, incompressible tori of minimal torus division number?

\end{question}

\noindent Our previous discussion of $T^2 \times [0,\infty)$, $T^2 \times \R$, and $S^1 \times \R^2$ proves
Theorem~\ref{thm: uncountableembeddings}.


\section{Proof of Theorem~\ref{thm: uncountable} and Threorem~\ref{thm: otuncountable}}

\begin{figure}

\includegraphics{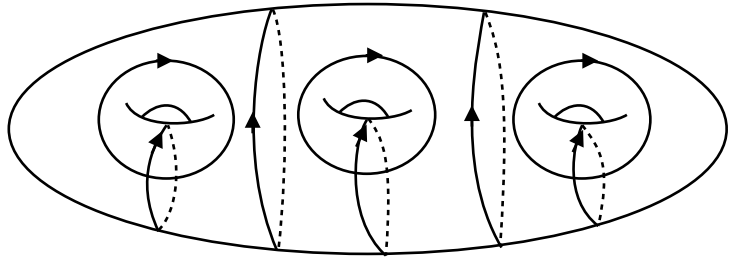}   

\caption{For $1 \leq i \leq 3$, let the $\alpha_i$ be the half-hidden, nonseparating, simple, closed curves and let $\beta_i$ be the nonseparating curves such that $\alpha_i \cdot \beta_i = 1$ (with subscript increasing from left to right).  Let the $\lambda_j$ be the two separating curves again labeled left to right.}\label{fig: lambda}
\end{figure}

Before beginning the proof of Theorem~\ref{thm: uncountable}, we prove a result which allows us to choose the dividing set on $\bdry M$ nicely.  Let $\Sigma$ be a genus $n$ surface.  In Figure~\ref{fig: lambda}, we specify $\alpha_i$, $\beta_i$, and $\lambda_j$ for a genus $3$ surface.  For a higher genus $\Sg$, make the analogous specification.

\begin{lemma} \label{lemma: kernel}

Let $M$ be any $3$-manifold with connected boundary of genus $n$.  Let $K$ be the kernel of the map $H_1(\partial M; \mathbb{Q}) \rightarrow
H_1(M; \mathbb{Q})$ induced from inclusion.  There exists an identification of $\partial M$ with $\Sg$ such that the $\alpha_i$ form a basis for $K \subset H_1(\partial M; \mathbb{Q})$ as vector space over $\mathbb{Q}$.  Moreover, there exist integers $n_i$ and embedded, orientable surfaces $\Sg_i$ such that $\partial \Sg_i$ consists of $n_i$ parallel copies of $\alpha_i$.\end{lemma}

\begin{proof}

Let $S_1$ be the first cutting surface in a Haken decomposition for $M$.  We may assume that no collection of components of $\partial S_1$ is separating in $\partial M$ and that $S_1$ is orientable ~\cite{He}.  We may also assume that $\partial S_1$ consists of parallel copies of a nonseparating, simple closed curve which we identify with $\alpha_1$.  If $\partial S_1$ is not all parallel, then choose two boundary components $b_1$ and $b_2$.  There exists an arc $\mu$ joining the $b_i$ which does not intersect any other components of $\partial S_1$.  Let $A$ be a small annular neighborhood of $\mu$.  Since $\partial S_1$ is nonseparating, we can choose $\mu$ so that $S_1 \cup A$ is an oriented surface with the $b_i$ replaced by a new boundary component homologous to $b_1 + b_2$.  We can continue this process until the boundary components of $S_1$ consist of $n_1$ copies of simple closed curve which we identify with $\alpha_1$.  Form a new $3$-manifold $M_1$ by attaching a $2$-handle $H_1$ to $\partial M$ along $\alpha_1$.  Let $S_2$ be the first surface in a Haken decomposition for $M_1$.  We may assume that $\partial S_2$ consists of $m_2$ copies of a nonseparating, simple closed curve $\gamma \subset \partial M_1$ which do not intersect the two disks $\partial H_1 \cap \partial M_1$.  Since $\partial S_2 \subset M$, we can identify $\gamma$ with $\alpha_2$.  Note that $S_2$ may intersect $H_1$.  If we cannot isotop the interior of $S_2$ to be disjoint from $H_1$, then we may assume that the intersection consists of $k$ disjoint disks $D_i$ on $S_2$.  Moreover, we can assume that the disks all have the same sign of intersection with the cocore of $H_1$.  For, if two disks had different signs of intersection, then we could find two adjacent such disks, remove the disks, and identify the boundaries to reduce the intersection of $S_2$ with $H_1$.  Note that $\partial S_1$ consists of $n_1$ copies of the attaching curve for $H_1$.  Therefore, we can take $k$ copies of $S_1$ and $n_1$ copies of $S_2$, remove the $kn_1$ disks $kS_2 \cap H_1$ from $n_1S_2$, and use the $kn_1$ boundary components of $kS_1$ to cap off these boundary components, possibly reversing the orientation of $S_1$ if necessary.  This operation shows that the class $n_2\alpha_2 \in K$, where $n_2 = n_1m_2$.  Attach another handle $H_2$ to $M_1$ along $\alpha_2$ to form a new manifold $M_2$.  Continuing in this fashion, we find $n$ integers $n_i$ and an identification of $\partial M$ with $\Sg$ such that $n_i\alpha_i \in K$.  The $\alpha_i$ are clearly linearly independent and thus generate $K$ since $dim_{\mathbb{Q}}(K) = n$ \cite{He}.\end{proof}

\begin{figure}

\includegraphics{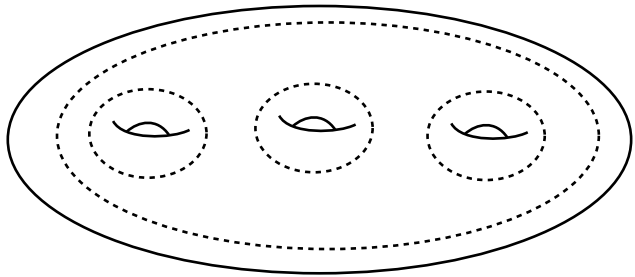}   

\caption{The collection of curves $\Gamma$ is diffeomorphic to the collection of curves shown above.}\label{fig: taut}

\end{figure}

Given any $3$-manifold with connected boundary, we identify $\partial M$ with the genus $n$ surface $\Sg$ as specified in Lemma~\ref{lemma: taut}.  We now describe the collection of curves $\Gamma \subset \partial M$ which will be the dividing set of a universally tight contact structure on $M$.  Let $\gamma_1$ be a simple closed curve homologous to $\alpha_1 -2\beta_1$ and let $\gamma_i$ be a simple closed curve homolgous to $\alpha_i -\beta_i$ for $2 \leq i \leq n$.  Finally, let $\gamma_{n+1}$ be a simple closed curve homologous to $-(\gamma_1 + \dots + \gamma_n)$.  Note that this collection of curves is diffeomorphic to the collection of curves shown in Figure~\ref{fig: taut}.

\begin{lemma} \label{lemma: taut}

Let $M$ be any irreducible $3$-manifold with connected boundary.  Then there exists a universally tight contact structure on $M$ such that $\partial M$ is convex and $\Gamma$ divides $\partial M$.  

\end{lemma}

\bp

Let $(M, \gamma)$ be the sutured $3$-manifold with annular sutures $s(\gamma) = \Gamma$.  We will show that $(M, \gamma)$
is a taut sutured $3$-manifold.  We then invoke the result in \cite{HKM2} which says that $M$ also supports a universally
tight contact structure with $\partial M$ convex and $\Gamma_{\partial M} = \Gamma$. 

To prove that $(M, \gamma)$ is taut, it suffices to show that $M$ is irreducible, $R(\gamma)$ is Thurston norm-minimizing
in $H_2(M, \gamma)$ among all other orientable surfaces in the same relative homology class, and $R(\gamma)$ is
incompressible in $M$.  By assumption, $M$ is irreducible.  We now show $R(\gamma)$ is incompressible.  Suppose not.  Then the Loop Theorem \cite{He} says that there exists an embedded disk $(D, \partial D) \subset (M, \partial M)$ such that $\partial D$ is homotopically nontrivial in $R(\gamma)$.  Since $R(\gamma)$ consists of two planar surfaces and $\partial D$ is embedded, $\partial D$ must also be homologous to $\pm (\gamma_{i_1} + \dots + \gamma_{i_j} )$ where $1 \leq i_1, i_j \leq n$ are distinct.  There exist $q_i \in \mathbb{Q}$ such that $\pm (\gamma_{i_1} + \dots + \gamma_{i_j} ) = q_1\alpha_1 + \dots + q_n\alpha_n$ since $\partial D$ is nulhomologous in $M$.  Take the intersection pairing of each side with $\alpha_i$ to arrive at a contradiction.  

\begin{figure}

\includegraphics{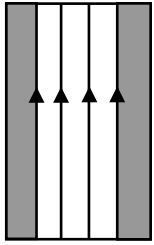}  

\caption{The white region is an annular suture.  The grey region is $R(\gamma)$.  The
two vertical lines in the annulus are boundary curves of $S$ with orientation induced from $S$. The arrows on $\partial
R(\gamma)$ denote the orientation induced from $R(\gamma)$.}\label{fig: orientations}

\end{figure}

We show that $R(\gamma)$ is Thurston norm-minimizing in $H_2(M, \gamma)$.  Let $S = \cup S_i$ be any orientable surface
homologous to $R(\gamma)$ in $H_2(M, \gamma)$.  Without loss of generality, we assume that $\partial S \subset int (A(\gamma))$.  Fix an annulus $A(s)
\subset A(\gamma)$ about the suture $s$ ($s$ is a homologically nontrivial simple closed curve in $A(s)$).  Note that
$\partial R(\gamma)$ intersects $A(s)$ in two oriented circles isotopic to $s$, where one comes from $R_{+}(\gamma)$ and
the other comes from $R_{-} (\gamma)$.  These circles must have the same orientation since the orientation of
$R_{+}(\gamma)$ agrees with the orientation on $\partial M$ and the orientation on $R_{-}(\gamma)$ does not.  Consider
the intersection of $S$ with $A(s)$.  If any two curves of $\partial S \cap A(s)$ have opposite orientation induced from
$S$, then we can find two such curves which are adjacent.  We then identify these curves and isotop them off of $\partial
M$ to reduce the number of boundary components of $S$.  We continue this procedure until $\partial S \cap A(s)$ consists of
two curves with the same orientation, which agrees with the orientation of $\partial A(s)$ induced from $R(\gamma)$.  Note
that the orientation on and number of these remaining curves in $\partial S \cap A(s)$ is completely determined by the
assumption that $[S] = [R(\gamma)]$ in $H_2(M, \gamma)$.  To summarize, we may assume that $\partial S$ intersects
each annulus of $A(\gamma)$ in exactly two essential curves with the same orientation induced from $S$, which agrees with
the orientation of the boundary of the annulus induced $R(\gamma)$ (see Figure~\ref{fig: orientations}).

We assume that our curves are exactly as in Figure~\ref{fig: taut}.  Recall that $\partial R(\gamma) = \cup_{j=1}^{n+1}\gamma_j \cup \cup_{j=1}^{n+1}\gamma_j $.  Let $S_i$ be a component of $S$.  We now show that $\partial S_i = \cup_{j=1}^{n+1}\gamma_j$ or $\partial
S_i = \cup_{j=1}^{n+1}\gamma_j \cup \cup_{j=1}^{n+1}\gamma_j$ as oriented manifolds.  Note that $\partial S_i$ is the union of some subset of the oriented curves $\{\gamma_1, \gamma_1, \dots, \gamma_{n+1},
\gamma_{n+1}\}$.  Since $\partial S_i \subset K$, $\partial S_i = q_1\alpha_1 + \dots + q_n\alpha_n$.  For $1 \leq j \leq n$, take the intersection pairing of both sides of this expression with $\alpha_j$ to see that $\gamma_j$ and $\gamma_{n+1}$ must occur together (if they occur at all) in $\partial S_i$.  This shows that $\partial S_i = \cup_{j=1}^{n+1}\gamma_j$ or $\partial
S_i = \cup_{j=1}^{n+1}\gamma_j \cup \cup_{j=1}^{n+1}\gamma_j$.  We say such surfaces are of type $I$ or $II$, respectively.  If $S_i$ is of type $I$, then $S$ consists of two such surfaces, and if $S_i$ is of type $II$, then $S=S_i$.  In either case, $x(S) \geq x(R(\gamma))$, with equality when $S$ is planar.\ep


\subsection{Construction of the contact structures when $\partial M$ is connected}\label{construction}


Let $(M, \eta)$ be the universally tight contact manifold given by Lemma~\ref{lemma: taut}.  When we refer to well-behaved surfaces, we will mean well-behaved with respect to $(M, \eta; \partial M, \lambda_1, \{\alpha_1, \beta_1\})$.  Let $S_1$ be the first cutting surface in a hierarchy for $M$ with boundary $\alpha_1$.  Recall that in the proof of Lemma~\ref{lemma: kernel} we chose $S_1$ so that $\partial S_1$ consists of $n_1$ copies of $\alpha_1$, so $S_1$ is well-groomed.  Via the correspondence between
sutured manifold decompositions and convex decompositions, we may assume that $S_1$ is the first cutting surface in a
convex decomposition for $M$ and has $\partial$-parallel dividing curves (see \cite{HKM2}).  Since $tb(S_1) \leq -2$, there is a bypass abutting $\partial M$ along $\alpha_1$.  
After attaching this bypass to $\partial M$, we have a $\Sg \times [0,1]$ slice
with convex boundary, where $\Sigma$ is a genus $n$ surface, $n$ is the genus of the boundary of $M$, and $\Sg \times
\{1\} = \bdry M$.  Let $(Y, \eta)$ denote this contact manifold.  Note that after attaching this bypass, the dividing
curves consist of $n$ $(-1,1)$ curves on each of the tori summands and another simple closed curve which is homologous to
the sum of the other $n$.

We now construct an embedding of $Y$ into $S^3$ with the standard tight contact structure.  Fix $g$ disjoint Darboux balls
in $S^3$ labeled $B_i$, where $g$ is the genus of the slice $Y$.  In $B_1$, we have a convex torus $T_1$ with slope $-2$. 
In each of the remaining $B_i$, we have a convex torus with slope $-1$.  On $T_1$, LeRP a curve $m_1$ which bounds a disk
in $T_1$ containing a single arc of the dividing set.  On each of the other $T_i$, LeRP a curve $l_i$ containing a disk in
$T_i$ with a single arc of $\Gamma_{T_i}$ and LeRP a curve $m_i$ which is disjoint from $l_i$ and bounds a disk with a
single arc of the same dividing curve that $l_i$ intersects. Now, remove the disks bounded by the $l_i$ and $m_i$ on $T_i$
and join $l_i$ to $m_{i+1}$ by a convex annulus $A_i$.  This yields a convex genus $n$ surface.  Inside $B_1$, we have a
compressing disk for $T_1$.  By the Imbalance principle, there is a bypass along this compressing disk.  Attaching this
bypass yields the desired embedding of $Y$.  Note that we can arrange for the sign of this bypass to agree with the sign of
the bypass we attached to $\bdry M$.

Fix a real number $r \in (-2,-1)$.  Let $q_i$ be an infinite sequence of rationals constructed in
Section~\ref{tighttoricends} such that $q_1 = -1$ and $q_i \neq r$.

\begin{prop} \label{prop: wellbehaved}

There exists a nested sequence $\Sg_i \subset Y = \Sg \times [0,1]$ of well-behaved surfaces such that
$\slope(\Sg_i) = q_i$ and $\Sg_1 = \Sg \times \{0\}$. 

\end{prop}

\begin{proof}

We will prove our results for the embedding of $Y \subset S^3$.  LeRP copies $l_i$ of $\lambda_1$ on $\Sg \times \{i\}$
such that $\tb(l_i) = -1$.  Let $A \subset Y$ be a convex annulus between $l_0$ and $l_1$.  $l_i$ separates $\Sg \times
\{i\}$ into a punctured torus $P_i$ and a punctured genus $n-1$ surface.  Cap off the $P_i$ in $S^3$ with convex disks
$D_i$ to obtain tori $T_i$ such that $\slope(T_1) = -2$ and $\slope(T_0) = -1$.  There exists an incompressible torus $T$
in the toric annulus bounded by the $T_i$ such that $\dv(T) = 1$ and $\slope(T) = q_2$ \cite{Ho2}.  Let $d_2$ be a
Legendrian divide on $T$.  $d_2$ can be Legendrian isotoped within the toric annulus bounded by the $T_i$ so that is
does not intersect the $D^2 \times [0,1]$ we used to cap off the thickened punctured torus bounded by $P_0 \cup P_1 \cup
A$.  This can be seen be working in a model for $D^2 \times [0,1]$, a standard neighborhood of a Legendrian arc. Hence, there exists a Legendrian isotopy taking $d_2$ to a curve in $Y$ that is homologous to $a_2\alpha_1 +
b_2\beta_1$, where $q_2 = b_2/a_2$.  LeRP a curve $d_2^{\prime}$ in the same homology class on $\Sg_1$ such that
$d_2^{\prime} \cap \Gamma_{\Sg_1}$ is minimal.  Let $A_2 \subset Y$ be a convex annulus between $d_2$ and
$d_2^{\prime}$.  By our choice of $d_2$, $\Gamma_{A_2} \cap d_2 = \emptyset$.  Attaching $A_1$ to $\Sg_1$
yields $\Sigma_2$.  Now, repeat the previous argument for $q_3$ and the slice bounded by $\Sg_2$ and $\Sg \times \{1\}$ to
obtain $\Sg_3$.  These surfaces are nested and well-behaved by construction.\end{proof}

Let $\Sg_i$ be as in Proposition~\ref{prop: wellbehaved}.  Let $(Y_i, \eta)$ be the genus $n$ slice bounded by $\Sg_i$ and
$\Sg_{i+1}$ in $Y$.  Construct a contact structure $\eta$ on $\Sg \times [0, \infty)$ by taking $\Sg \times [i, i+1]$ to be
$Y_i$.  Let $(V, \eta_r)$ be obtained from $(M, \eta)$ by peeling off $Y \setminus \Sg_1$ from $(M, \eta)$ and attaching
$(\Sg \times [0, \infty), \eta)$ in the obvious way.  Note that $(V, \eta_r)$ is tight by construction and embeds into
$(M, \eta)$.  

\begin{lemma}\label{lem: isotopy} 

Let $s,t \in (-2,-1)$.  Then $(V, \eta_s)$ and $(V, \eta_t)$ are in the same isotopy class of contact structures.

\end{lemma}

\begin{proof}
 
There exists a convex surface $S \subset V$ such that $V \setminus S = V' \cup S \times (0,\infty)$, where $V'$ is diffeomorphic
to $V$ and $\eta_s|_{V' \cup S} = \eta_t|_{V' \cup S}$.  This follows from the construction of $(V, \eta_s)$ and $(V, \eta_t)$.  
We claim that $\eta_s|_{S \times [0,\infty)}$ and $\eta_t|_{S \times [0,\infty)}$ are isotopic rel $S \times 0$.  We can assume that $S \times [0,1)$ is a one-sided vertically invariant neighborhood of our convex surface $S \times 0$.  Hence, in particular, $\eta_t$ and $\eta_s$ agree on $S \times [0,1)$.  Form a new contact structure $\eta_t^\lambda$ as follows:  Extend the vertically invariant neighborhood $S \times [0,1)$ of $\eta_t$ to $S \times [0, \lambda)$, and on $S \times [\lambda, \infty)$ take $\eta_t^\lambda$ to be $\eta_t |_{S \times [1,\infty)}$.  Define $\eta_s^\lambda$ similarly.  By construction, $\eta_t^1 = \eta_s^1$.  Hence, $(V, \eta_s)$ and $(V, \eta_t)$ are in the same isotopy class of contact structures.\end{proof}


\subsection{Proof of Theorem~\ref{thm: uncountable} and Theorem~\ref{thm: otuncountable} when $\partial M$ is connected}\label{incomp}


In order to show that $V$ supports uncountably many tight contact structures that are not contactomorphic, we will first
show that the $(V, \eta_s)$ are distinct up to proper isotopy.  Theorem~\ref{thm: uncountable} then follows immediately since the
mapping class group of any $3$-manifold with boundary is countable (\cite{McC}).  To achieve this, we use the idea of the
slope at infinity introduced in Section~\ref{invariant}.  

\begin{prop} \label{prop: beta}

The net $\slope \co \mathcal{C}(Ends(V, \eta_s; \partial M)) \rightarrow \R \cup \{\infty\}$ is convergent, so the slope at
infinity is defined.  Moreover, the slope at infinity is $s$ for all $s \in (-2,-1)$.

\end{prop}

\begin{proof}

We first show that there is an $E \in Ends(V, \eta_s)$ such that for all $F \subset E$, $\slope(F) \leq s$.  Choose $E
\subset \intr(Y)$.  We will be now working in $S^3$.  Let $F \subset E$ and suppose for contradiction that $\slope(F) > s$.  Then, there
exists $\Sigma \in \mathcal{C}(E)$ such that $\slope(\Sigma) > s$.  Let $\Sg_i$ be the family of surfaces given by Proposition~\ref{prop: wellbehaved}.  There exists an $i$ such that $\Sg$ is contained in the
genus $n$ slice bounded by $\Sg_1$ and $\Sg_i$.  LeRP a copy of $\lambda_1$ on $\Sg$, $\Sg_1$ and $\Sg_i$ and cap off the
punctured tori bounded by these curves with convex disks.  This yields a toric annulus $T^2 \times [0,1] \subset S^3$ 
which contains a convex, incompressible torus $T$ such that $\slope(T) > \slope(T^2 \times \{1\})$.  No such $T^2 \times I$
can exist in $S^3$ (see \cite{Ho2}).  Therefore, such a $\Sigma$ could not exist.  Similarly, one can show that $\slope(F) < s$ leads to a contradiction.  The existence of the family $\Sg_i$ now implies that the slope at infinity is $s$.\end{proof}

By the proper isotopy invariance of the slope at infinity, there are uncountably many tight contact structures that are not
properly isotopic on $V$.  This concludes the proof of Theorem~\ref{thm: uncountable} in the case of an irreducible $M$
with connected boundary.  The proof of Theorem~\ref{thm: otuncountable} is now immediate.  For each $\eta_s$, simply choose a transverse curve in
$V$ and introduce a Lutz twist.  Since the contact structures is identical outside of a compact set, the slope at infinity is unchanged.


\subsection{Proof of Theorem~\ref{thm: uncountable} when $\partial M$ is disconnected}\label{gencase}


Before proceeding with the proof, we will need the following technical result.

\begin{lemma}\label{lemma: fillings}
For every integer $n$, there exists an irreducible $3$-manifold $M_n$ with connected, incompressible boundary of genus $n$.
\end{lemma}

\bp
Let $\Sg_g$ be an orientable surface of genus $g$.  If $n = 2m$, let $F \subset \Sg_n$ be a once punctured genus $m$ surface.  Form a manifold $M_n$ by identifying $F \times \{0\}$ with $F \times \{1\}$ on $\Sg_n \times [0,1]$.  It is straightforward to show that $\Sg_n \times \{0\}$ and $\Sg_n \times \{1\}$ are incompressible in $M_n$.  Using the incompressibility of these surfaces and the irreducibility of $\Sg_n \times [0,1]$, it is routine to show that $M_n$ has incompressible boundary and is irreducible.  If $n = 2m-1$, let $F \subset \Sg_m$ be an annular neighborhood of a nonseparating simple closed curve.  Form a manifold $M_n$ by identifying $F \times \{0\}$ with $F \times \{1\}$ on $\Sg_m \times [0,1]$.  It is again straightforward to show that $\Sg_m \times \{0\}$ and $\Sg_n \times \{1\}$ are incompressible in $M_n$.  Irreducibility and incompressibility of the boundary follow as before.\ep

Let $\partial M = \cup_{i=1}^n S_i$ where the $S_i$ are the connected components of $\partial M$ and $S_1$ is of nonzero genus.  Let $S_j$ be any component different from $S_1$.  If $S_j$ is compressible, compress it, and continue doing so until we have a collection of spheres and incompressible surfaces.  We are now in the situation where every boundary component, besides possibly $S_1$, is incompressible or a sphere.  Fill in each sphere with a ball and onto each incompressible component of genus $n$, excluding $S_1$ if it happens to be incompressible, glue in an irreducible manifold with connected, incompressible boundary of genus $n$ (such manifolds exist by Lemma~\ref{lemma: fillings}).  It is straightforward to show that the resulting manifold is irreducible since we are gluing irreducible manifolds along incompressible surfaces.  Call the resulting manifold $M'$.  We are now in the case of connected boundary.  Put a tight contact structure on $M'$ as before and attach a bypass along $\alpha_1$ so that we have factored off a $\partial M' \times [0,1]$ slice $Y$.  Topologically, the closure of $M' \setminus Y$ is again $M'$.  Without intersecting $Y$, remove each of the manifolds we glued in after perturbing the gluing surfaces to be convex.  Reconstruct $M$ by gluing the boundary components back together along the compressing disks.  To ensure that the resulting manifold is tight, choose the compressing disks to be convex with Legendrian boundary and with a single arc in the dividing set \cite{Ho1}.  We now have a tight contact structure on $M$ and a bypass layer $Y$ along $S_1$ which is identical to the case of connected boundary.  To form the $(V, \eta_s)$ remove all the boundary components except for $S_1$ and construct the ends in $Y$ as before.  The calculation of the slope at infinity is identical to the case of connected boundary.  As in the case of connected boundary, the proof of Theorem~\ref{thm: otuncountable} is immediate after introducing a Lutz twist along a transverse curve in
$V$.  

\s\s\n  {\em Acknowledgements.}  

Thank you to John Etnyre, my advisor, and Stephan Sch\"{o}nenberger for helpful conversations and for reading drafts of this work.  Also, thank you to Ko Honda, Will Kazez, and Gordana Mati\'c for their comments and questions during my talk at the Georgia Topology Conference and the conversations that followed.

\end{document}